\documentclass{amsart}
\usepackage{amssymb}

\newtheorem{theorem}{Theorem}
\newtheorem{proposition}[theorem]{Proposition}
\newtheorem{lemma}[theorem]{Lemma}
\newtheorem{remark}[theorem]{Remark}

\newcommand{\R}{\mathbb{R}}
\newcommand{\Z}{\mathbb{Z}}
\newcommand{\A}{\mathcal{A}}
\newcommand{\B}{\mathcal{B}}

\newcommand{\form}{\overline{a}}
\newcommand{\Schaefer}{\mathcal{S}}
\newcommand{\abs}[1]{\left|#1\right|}
\newcommand{\norm}[1]{\left\| #1 \right\|}
\newcommand{\skprd}[2]{\left\langle #1 , #2\right\rangle}

\DeclareMathOperator{\relCap}{Cap_{\overline{\Omega}}}
\DeclareMathOperator{\Lmax}{L_{max}(\varepsilon,\lambda_1)}

\newcommand{\vup}{\overline{v}}
\newcommand{\vdown}{\underline{v}}

\begin{document}

\title[Semilinear elliptic equations on arbitrary open sets]{Existence of solutions for semilinear elliptic boundary value problems on arbitrary open sets}
\author[Reinhard Stahn]{Reinhard Stahn *}

\begin{abstract}
We show the existence of a weak solution of a semilinear elliptic Dirichlet problem on an arbitrary open set $\Omega$. We make no assumptions about the open set $\Omega$ and very mild regularity assumptions on the semilinearity $f$, plus a coerciveness assumption which depends on the optimal Poincar\'{e}-Steklov constant $\lambda_1$. The proof is based on Schaefer's fixed point theorem applied to a sequence of truncated problems. We state a simple uniqueness result. We also generalize the results to Robin boundary conditions. 
\end{abstract}

\maketitle

{\raggedright
}

{\let\thefootnote\relax\footnotetext{*Fachrichtung Mathematik, Institut f\"{u}r Analysis, Technische Universit\"{a}t Dresden, 01062, Dresden, Germany. Email: \textit{Reinhard.Stahn@mailbox.tu-dresden.de}}}


\section{Introduction}\label{Einleitung}
The question of existence of (weak) solutions for a boundary value problem like
\begin{align}
  \label{General BVP}
  \left\{
  \begin{array}{rlr}
    -\Delta u(x) & = f(x,u(x),\nabla u(x)) & (x\in\Omega)\\
    \B u(x) & = 0 & (x\in \partial \Omega).
  \end{array} 
	\right.
\end{align}
is a classical problem. Here $\B u(x)=0 \, (x\in \partial \Omega)$ is simply an abbreviation for specific boundary conditions (eg. Dirichlet, Neumann or Robin boundary conditions). However, in older publications (e.g. \cite{AmannCrandall1978,KazdanKramer1978,Pokhozhaev1980} and references therein) it seems to be standard to assume that $\Omega$ is bounded and has sufficiently smooth boundary. On the other hand for a long time it is well known that Dirichlet (and Neumann) boundary conditions can easily be defined on arbitrary open sets in a weak sense. But until now a general existence theorem for weak solutions of (\ref{General BVP}) on arbitrary open sets (that means in particular possibly unbounded, not smooth, not connected), possibly irregular semilinearities $f$ (that means in particular possibly not Lipschitz continuous in $u$ or $\nabla u$) and general boundary conditions seems to be unknown. 

In this paper we state and prove such a theorem with the help of a fixed point theorem in a locally convex space. For simplicity we restrict ourselves to the case of Dirichlet boundary condition. A short discussion on how the results extend to other boundary conditions is included at the end of the paper. To find solutions of (\ref{General BVP}) with the help of fixed point methods is standard. However, in the majority of textbooks a Banach space setting is presented (cf. \cite{DrabekMilota2007,Evans2010,GilbargTrudinger2001}).


\section{Assumptions and main results}\label{Introduction}
Let $\Omega\subseteq\R^d$ be an arbitrary open set. We consider the following Dirichlet problem:

\begin{align}
  \label{DP}
  \left\{
  \begin{array}{rlr}
    -\Delta u(x) & = f(x,u(x),\nabla u(x)) & (x\in\Omega)\\
    u(x) & = 0 & (x\in \partial \Omega).
  \end{array} 
	\right.
\end{align}
We call $u$ a \textit{weak solution} of this problem iff $u\in H^1_0(\Omega)$ such that $f(x,u,\nabla u)\in L_2(\Omega)$ and 
\begin{align*}
	\int_{\Omega} \nabla u\cdot \nabla \varphi dx = \int_{\Omega} f(x,u,\nabla u)\varphi dx \quad \forall \varphi \in C_c^{\infty}(\Omega).
\end{align*}
In the following we show under which assumptions on the semilinearity $f$ we can prove a priori bounds, existence and uniqueness of weak solutions. 

\subsection{Assumptions and notation}\label{Assumptions}
We assume that $f:\Omega\times\R\times\R^d\rightarrow\R$ is a \textit{Caratheodory function}. This means that $f=f(x,s,\xi)$ is measurable as a function in $x$ when $s$ and $\xi$ are fixed, and is jointly continuous as a function in $s$ and $\xi$ when $x$ is fixed, for almost all $x$.
The function $f$ should satisfy a \textit{coerciveness} and a \textit{growth condition}
\begin{align}
  \label{coercive} 
	f(x,s,\xi)s &\leq (\lambda_1-\varepsilon)s^2 + L\abs{\xi s} + h(x)\abs{s} \quad \forall x\in\Omega, s\in\R, \xi\in\R^d , \\
	\label{Wachstumsbedingung}
	f(x,s,\xi)s &\geq -\gamma(\abs{s})\abs{s} - L_0\abs{\xi s} - h_0(x)\abs{s} \quad \forall x\in\Omega, s\in\R, \xi\in\R^d .
\end{align}
Here $\lambda_1$ is the optimal Poincar\'{e}-Steklov constant for the Dirichlet-Laplace operator, i.e. $\lambda_1\geq 0$ is the largest real number such that $\lambda_1\norm{u}^2_{L_2(\Omega)}\leq\norm{\nabla u}^2_{L_2(\Omega)}$ is true for all $u\in H^1_0(\Omega)$. The positive constant $L$ has to satisfy 
\begin{align}\label{L Bedingung}
	L < \Lmax :=
  \begin{cases}
   \varepsilon/\sqrt{\lambda_1} & \text{if } \varepsilon \leq 2\lambda_1, \\
   2\sqrt{\varepsilon-\lambda_1} & \text{if } \varepsilon \geq 2\lambda_1.
  \end{cases} 
\end{align}
Furthermore $\varepsilon>0$ and $h\geq0$ with $h\in L_2\cap L_q(\Omega)$ for some $q\geq 2$ which (for simplicity) is not equal to $d/2$. The number $q$ will serve as a parameter. $L_0\geq L$ is an arbitrary constant and $h_0\in L_2(\Omega)$ with $h_0\geq0$. The monotone increasing function $\gamma:[0,\infty)\rightarrow [0,\infty)$ is assumed to satisfy
\begin{align}\label{Gamma bei 0}
	\limsup_{s\rightarrow 0} \frac{\gamma(s)}{s} < \infty,
\end{align}
and in case of $q< d/2$ also
\begin{align}\label{Gamma bei unendlich}
	\limsup_{s\rightarrow +\infty} \frac{\gamma(s)}{s^{q^{**}/2}} < \infty
\end{align}
for the positive real number $q^{**} = qd/(d-2q)$. For future use we also define $q^{**}=\infty$ if $q>d/2$. The bigger the parameter $q$, the more restrictive is the coerciveness condition (\ref{coercive}) but the less restrictive is the growth condition (\ref{Wachstumsbedingung}). Note that the condition (\ref{Gamma bei 0}) is not needed if the measure of $\Omega$ is finite. In the case when (\ref{Gamma bei 0}) is not satisfied but the measure of $\Omega$ is finite we could change $\gamma(s)$ to $(\gamma(s)-\gamma(1))_+$ and add the additional constant $\gamma(1)$ to the function $h_0$ without touching its $L_2$-integrability. Here and in the following for a real number $a$ we define $a_+:=\max\{a,0\}$.

We will also consider the following two Dirichlet problems on $\Omega$
\begin{align}
  \label{DP0}
  \left\{
  \begin{array}{rlr}
    -\Delta v & = (\lambda_1-\varepsilon) v - L\abs{\nabla v} - h(x) & (x\in\Omega)\\
    v(x) & = 0 & (x\in\partial \Omega)
  \end{array} 
	\right.
\end{align}
and
\begin{align}
  \label{DP1}
  \left\{
  \begin{array}{rlr}
    -\Delta v & = (\lambda_1-\varepsilon) v + L\abs{\nabla v} + h(x) & (x\in\Omega)\\
    v(x) & = 0 & (x\in\partial \Omega).
  \end{array} 
	\right.
\end{align}
In Section \ref{Domination in Lr} we prove that solutions of these two equations exist and are unique. Therefore by $\vdown, \vup$ let us denote the solutions of (\ref{DP0}) and (\ref{DP1}), respectively. We will also see that $\vdown \leq 0 \leq \vup$. By (\ref{coercive}) this implies $-\Delta\vdown\leq f(x,\vdown,\nabla\vdown)$ and $-\Delta\vup\geq f(x,\vup,\nabla\vup)$ that is, these two functions are sub- and supersolutions of (\ref{DP}).

\subsection{Main results}\label{Main results}
The first step to prove existence of a weak solution of (\ref{DP}) is to prove a priori estimates for hypothetical solutions of a class of semilinear problems: Let $\mu\geq \max\{0,\varepsilon-\lambda_1\}$, $0<t\leq 1$ and $\omega$ be an open subset of $\Omega$. Moreover, let $v_0\leq0\leq v_1$ be two measurable functions on $\Omega$. Consider
\begin{align}
  \label{DPklasse}
  \left\{
  \begin{array}{rlr}
    -\Delta u + \mu u & = t\underbrace{(f(x,\sigma(x,u),\nabla u)+\mu \sigma(x,u))}_{=:b_{\sigma}(x,u,\nabla u)}\chi_{\omega}(x) & (x\in\Omega)\\
    u(x) & = 0 & (x\in\partial \Omega)
  \end{array} 
	\right.
\end{align}
where $\sigma(x,s)=\max\{v_0(x),\min\{s,v_1(x)\}\}$. 

\begin{theorem}[a priori bounds in $H^1_0$]\label{H1}
Assume that $f$ satisfies (only) the coerciveness condition (\ref{coercive}). Let $u$ be a weak solution of (\ref{DPklasse}). Then $u$ satisfies an a priori estimate in $H^1_0(\Omega)$. More precisely, there exists a constant $C>0$ such that
\begin{align*}
	\norm{u}_{H^1_0(\Omega)} \leq C\norm{h}_{L_2(\Omega)}.
\end{align*}
The constant depends only on $\lambda_1,\varepsilon$ and $(\Lmax-L)^{-1}$.
\end{theorem}
This theorem is the central argument in the proof of the existence of a weak solution of (\ref{DP}) in Section \ref{Proof of Theorem}.

\begin{theorem}[a priori bounds in $L_{q^{**}}$]\label{Lr}
Assume that $f$ satisfies (only) the condition (\ref{coercive}), where the term $\lambda_1-\varepsilon$ is replaced by any positive number $\delta$. Let $u$ be a weak solution of (\ref{DPklasse}). Then $u$ satisfies an a priori estimate in $L_{q^{**}}(\Omega)$. More precisely, there exists a constant $C>0$ such that
\begin{align*}
	\norm{u}_{L_{q^{**}}(\Omega)} \leq C\left(\norm{u}_{L_2(\Omega)} + \norm{h}_{L_q(\Omega)}\right).
\end{align*}
The constant depends only on $d,q,\delta$ and $L$.
\end{theorem}

Much more interesting than merely a bound for the $L_{q^{**}}$-norm is the fact that weak solutions are even a priori dominated by functions in the $L_{q^{**}}(\Omega)$. 

\begin{theorem}[a priori domination in $L_{q^{**}}$]\label{Dominationtheorem}
The Dirichlet problems (\ref{DP0}) and (\ref{DP1}) have unique weak solutions $\vdown\leq0$ and $\vup\geq0$ respectively. Moreover,
if we assume that $f$ satisfies (only) the coerciveness condition (\ref{coercive}) then
\begin{align*}
	\vdown \leq u\leq \vup 
\end{align*}
for every weak solution $u$ of (\ref{DPklasse}).
\end{theorem}
 
The main result of this paper is:

\begin{theorem}[existence]\label{Haupttheorem}
Under the assumptions in Section \ref{Assumptions} the Dirichlet problem (\ref{DP}) admits a weak solution. 
\end{theorem}

Note that only in this theorem, but not in the preceding three theorems we assume the validity of the growth condition (\ref{Wachstumsbedingung}). The assumption (\ref{L Bedingung}) on $L$ in (\ref{coercive}) might look a bit strange. But for our proof of the a priori bounds in $H^1_0$ it is needed. For a discussion of this assumption see section \ref{Gegenbeispiel}. It is not the main concern of this paper but whenever one can show an existence theorem the question of \textit{uniqueness} arises. In section \ref{Uniqueness} a simple condition is given when uniqueness holds. Section \ref{Robin BVPs} is devoted to Robin boundary conditions. 


\section{Schaefer's fixed point theorem}\label{Schaefer}
The proof of Theorem \ref{Haupttheorem} is mainly based on Schaefer's fixed point theorem (cf. \cite{Schaefer1955}). We do not apply the original formulation of Schaefer's theorem from his paper, but a slightly different version. Actually it does not matter which version we apply but we prefer the version presented below.
\begin{theorem}[Schaefer's fixed point theorem]\label{Schaefers Theorem}
Let $X$ be a locally convex Hausdorff space and $T:X\rightarrow X$ a continuous mapping. Let 
\begin{align*}
\Schaefer	= \{u\in X| \exists 0<t\leq 1: u=tTu \}.
\end{align*}
Let $\norm{\cdot}$ be a continuous semi-norm on $X$, let $\rho>0$ be a number and $K_{\rho}=\{u\in X: \norm{u}<\rho\}$. If 
\begin{itemize}
	\item[(i)] $\Schaefer \subseteq K_{\rho}$ and
	\item[(ii)] $TK_{\rho} \subseteq X$ is relatively compact,
\end{itemize}
then $T$ has a fixed point $u=Tu$.
\end{theorem}
A proof which is based on the Tychonoff-Schauder fixed point theorem (cf. \cite{Tikhonov1935}) can be found in \cite[Section 2]{ArendtChill2010}. Actually we could replace (i) by the more general condition $u\in\Schaefer\Rightarrow \norm{u}\neq\rho$. This shows the connection of Schaefer's theorem to the well known degree theory of Leray-Schauder.

\section{A priori bounds and domination}\label{Bounds and Domination}
We consider the class of boundary value problems (\ref{DPklasse}) and prove the Theorems \ref{H1}, \ref{Lr} and \ref{Dominationtheorem}. It is easy to generalize these Theorems also to sub- and supersolutions.

\subsection{A priori bounds in $H^1_0(\Omega)$}
\textit{Proof of Theorem \ref{H1}.} Let $u$ be a weak solution of (\ref{DPklasse}). We choose $u\in H^1_0(\Omega)$ as a test function for (\ref{DPklasse}). Since $\lambda_1+\mu-\varepsilon\geq0, \abs{s}\geq\abs{\sigma(x,s)}$ and the fact that $\sigma(x,s)$ and $s$ have the same sign for every $s\in\R$ and almost all $x\in\Omega$ by the coerciveness condition (\ref{coercive}) we may deduce
\begin{align} \label{Hnorm}
	\norm{\nabla u}_{L_2(\Omega)}^2 + \mu \norm{u}_{L_2(\Omega)}^2 \leq & t(\lambda_1 + \mu - \varepsilon)\norm{u}_{L_2(\omega)}^2 \\ \nonumber
	+ &tL\norm{\nabla u}_{L_2(\omega)}\norm{u}_{L_2(\omega)} \\ \nonumber
	+ &t\norm{h}_{L_2(\Omega)}\norm{u}_{L_2(\omega)}.
\end{align}
Let $\varepsilon_1,\varepsilon_2\geq0$ with $\varepsilon=\varepsilon_1+\varepsilon_2$ for which there is a $\delta_1\in(0,1]$ such that $\varepsilon_1=\delta_1 \lambda_1$. If $\lambda_1=0$ we choose $\delta_1=1$. Appropriately inserting the Poincar\'e-Steklov inequality into (\ref{Hnorm}) yields
\begin{align*}
	&((1-\delta_1)\lambda_1+\mu)\norm{u}_{L_2(\omega)}^2 + \delta_1\norm{\nabla u}_{L_2(\omega)}^2 \\
	&\leq t((1-\delta_1)\lambda_1+\mu-\varepsilon_2)\norm{u}_{L_2(\omega)}^2 \\
	&+ tL\norm{\nabla u}_{L_2(\omega)}\norm{u}_{L_2(\omega)} + t\norm{h}_{L_2(\Omega)}\norm{u}_{L_2(\omega)} .
\end{align*}
After rearranging the terms and dividing by $t$ we get
\begin{align*}
  &((t^{-1}-1)(\lambda_1-\varepsilon_1+\mu)+\varepsilon_2)\norm{u}_{L_2(\omega)}^2 + t^{-1}\delta_1\norm{\nabla u}_{L_2(\omega)}^2 \\
	&\leq L\norm{\nabla u}_{L_2(\omega)}\norm{u}_{L_2(\omega)} + \norm{h}_{L_2(\Omega)}\norm{u}_{L_2(\omega)} .
\end{align*}
Since $(t^{-1}-1)(\lambda_1-\varepsilon_1+\mu)\geq0$ and $t^{-1}\delta_1\geq\delta_1$ we deduce with the help of the inequality of arithmetic and geometric mean
\begin{align}\label{Bumm}
	&\varepsilon_2\norm{u}_{L_2(\omega)}^2 + \delta_1\norm{\nabla u}_{L_2(\omega)}^2 \\ \nonumber
	&\leq \delta_1\norm{\nabla u}_{L_2(\omega)}^2 + \frac{L^2}{4\delta_1}\norm{u}_{L_2(\omega)}^2 + \norm{h}_{L_2(\Omega)}\norm{u}_{L_2(\omega)}
\end{align}
i.e.
\begin{align*}
	4\delta_1\varepsilon_2 \norm{u}_{L_2(\omega)} \leq L^2\norm{u}_{L_2(\omega)} + 4\delta_1\norm{h}_{L_2(\Omega)}.
\end{align*}
If we set $\varepsilon_1=\min\{\varepsilon/2, \lambda_1\}$ we easily see that $4\delta_1\varepsilon_2>L^2$ and thus we deduce that there exists a $\rho_0>0$ such that $\norm{u}_{L_2(\omega)} \leq \rho_0 \norm{h}_{L_2(\Omega)}$. Inserting this in (\ref{Hnorm}) we get the desired a priori bound.
\hfill $\Box$

\subsection{Uniqueness.}\label{Uniqueness} A second view on the proof of the Theorem \ref{H1} gives rise to a uniqueness result for (\ref{DP}) if we assume the \textit{monotonicity} condition 
\begin{align}\label{Monotonie}
	(f(x,s_2,\xi_2)-f(x,s_1,\xi_1))(s_2-s_1) &\leq (\lambda_1-\varepsilon)(s_2-s_1)^2 + L\abs{\xi_2-\xi_1}\abs{s_2-s_1}.
\end{align}

\begin{proposition}\label{uniqueness Proposition}
Let $f$ satisfy the monotonicity condition (\ref{Monotonie}). Then the Dirichlet problem (\ref{DP}) has at most one weak solution.
\end{proposition}
\textit{Proof.} Let us assume that there exist two solutions $u_1$ and $u_2$. Then $u=u_2-u_1\in H^1_0(\Omega)$ can serve as a test function for (\ref{DP}) with $u$ replaced by $u_1$ or $u_2$. Then (\ref{Monotonie}) leads to 
\begin{align*}
	\norm{\nabla u}_{L_2(\Omega)}^2 
	&\leq (\lambda_1 - \varepsilon)\norm{u}_{L_2(\Omega)}^2 + L\norm{\nabla u}_{L_2(\Omega)}\norm{u}_{L_2(\Omega)} \\
	&\stackrel{(*)}{\leq} \norm{\nabla u}_{L_2(\Omega)}^2 + (L-\Lmax)\norm{\nabla u}_{L_2(\Omega)}\norm{u}_{L_2(\Omega)} .
\end{align*}
We deduce $u=0$. The second inequality $(*)$ follows from the Poincar\'{e}-Steklov inequality and 
\begin{align*}
	\Lmax\norm{\nabla u}_{L_2(\Omega)}\norm{u}_{L_2(\Omega)} 
	\leq \delta_1\norm{\nabla u}_{L_2(\omega)}^2 + \frac{\Lmax^2}{4\delta_1}\norm{u}_{L_2(\omega)}^2
\end{align*}
for $\delta_1=\varepsilon/(2\lambda_1)$ if $\varepsilon\leq 2\lambda_1$ and $\delta_1=1$ else.
\hfill $\Box$

\subsection{A priori bounds in $L_{q^{**}}(\Omega)$}
\textit{Proof of Theorem \ref{Lr}.} Theorem \ref{H1} already implies an a priori bound for the $L_{2^*}(\Omega)$-norm by the Sobolev embedding theorem. With the help of Moser's iteration method we also achieve bounds with respect to higher order Lebesgue norms. We learned the Moser iteration technique from the proof of \cite[Theorem 8.15]{GilbargTrudinger2001} and could generalize it to our situation. For equations of the form $-\Delta u=g\in L_q(\Omega)$ on bounded domains and Robin boundary conditions see also \cite[Theorem 4.1]{Daners2000}.

\textbf{(i) Making Moser iteration possible.} 
For $2\leq r \leq \infty$ and $2\leq p\leq q$ we define
\begin{align*}
	M(r) &= \max\{\norm{u}_{L_2(\Omega)}, \norm{u}_{L_r(\Omega)}\} \text{ and} \\
	\rho(p) &= M(p)+\norm{h}_{L_p(\Omega)}.
\end{align*}
By Lyapunov's interpolation inequality\footnote{Lyapunov's inequality \cite[Section 2.9]{HardyLittlewoodPolya1934}: $\norm{v}_{L_r(\Omega)}\leq\norm{v}_{L_p(\Omega)}^{1-\theta}\norm{v}_{L_q(\Omega)}^{\theta}$ for all measurable functions $v$ if $1/r=(1-\theta)/p+\theta/q$ where $0< r,p,q\leq\infty$ and $0\leq\theta\leq 1$.} for Lebesgue spaces we see that $M$ is an increasing function. Moreover $M:[2,\infty]\rightarrow[0,\infty]$ is continuous. A standard test function argument shows that the coerciveness condition (\ref{coercive}) implies
\begin{align}\label{Rekursionsungleichung}
  M(2^*\beta)^{2\beta} \leq C_1\beta^2 \rho(p) M(p'(2\beta-1))^{2\beta-1}
\end{align}
for all $\beta\geq 1$ and $2\leq p\leq q$. Here $2^*$ is equal to $2d/(d-2)$ if $d>2$, or some sufficiently big number else. The constant $C_1\geq 1$ only depends on $d,\delta$ and $L$. 

For the proof of (\ref{Rekursionsungleichung}) one can make the same ansatz as in the proof of \cite[Theorem 8.15]{GilbargTrudinger2001} with $k=0$. That is, we use a test function which is proportional to $\abs{u}^{2\beta-2}u$ for small values $\abs{u(x)}$ and proportional to $u$ for big values of $\abs{u(x)}$. The crucial point in the proof is the validity of the Sobolev embedding $H^1_0(\Omega)\hookrightarrow L_{2^*}(\Omega)$. 

\textbf{(ii) Moser iteration.} We distinguish the two cases $2q<d$ and $2q>d$.

\textit{Case $2q<d$:} For the beginning let us assume that $2\leq p< \max\{2^*,q\}$. By Sobolev's embedding theorem $\rho(p)<\infty$. Note that $2^*\beta=p'(2\beta-1)$ is equivalent to $2^*\beta=p^{**}$. Thus (\ref{Rekursionsungleichung}) implies $M(p^{**})<\infty$ and therefore
\begin{align}\label{ok}
	M(p^{**}) \leq C_1\left(\frac{p^{**}}{2^*}\right)^2 \left(M(p)+\norm{h}_{L_p(\Omega)}\right).
\end{align}
This yields $\rho(p^{**})<\infty$. By iterating the preceding argument we get $M(p^{**})<\infty$ for all $2\leq p<q$ and a \textit{uniform} estimate holds. By the continuity of $M$ this is also true for $p=q$. By Lyapunov's inequality at the cost of a bigger constant we can replace $M(q)$ on the right hand side of (\ref{ok}) by $\norm{u}_{L_2(\Omega)}$ and the claim follows. 

\textit{Case $2q>d$:} From the first case we already know that $M(r)<\infty$ for all $2\leq r<\infty$. Thus $\rho(q)<\infty$. Let us recursively define the increasing sequence $(\beta_n)$ by $\beta_0=1, \beta_{n}=\frac{1}{2}+\chi \beta_{n-1}$ where $\chi=\frac{2^*}{2q'}>1$. Then we deduce from (\ref{Rekursionsungleichung}) with $p=q$ that
\begin{align*}
  M(2^*\beta_{n}) \leq \beta_n^{\frac{1}{\beta_n}} (C_1\rho(q))^{\frac{1}{2\beta_n}} M(2^*\beta_{n-1})^{1-\frac{1}{2\beta_n}}
\end{align*}
Observe that $1-\frac{1}{2\beta_n}=\frac{\chi \beta_{n-1}}{\beta_n}$. This gives immediately
\begin{align*}
  M(2^*\beta_N) \leq \underbrace{\left(\prod_{m=1}^{N}\beta_m^{\chi^{N-m}}\right)^{\frac{1}{\beta_N}}}_{C(N,d,q)} (C_1\rho(q))^{1-\frac{\chi^N}{\beta_N}} M(2^*)^{\frac{\chi^N}{\beta_N}}.
\end{align*}
It is a simple exercise to show that
\begin{itemize}
	\item [(a)] $\chi^N/\beta_N$ converges from above to $\theta=(2\chi-2)/(2\chi-1)\in(0,1)$ and 
	\item [(b)] there exists a constant $C(d,q)>0$ such that $C(N,d,q)\leq C(d,q)$.
\end{itemize}
Thus
\begin{align*}
  M(2^*\beta_N) \leq C_2 \rho(q)^{1-\frac{\chi^N}{\beta_N}} M(2^*)^{\frac{\chi^N}{\beta_N}}.
\end{align*}
The constant $C_2$ depends on $d,\delta,L$ and $q$. Since $M$ is continuous we may let $N$ tend to infinity and deduce from (b) that
\begin{align*}
  M(\infty) \leq C_2 \rho(q)^{1-\theta} M(2^*)^{\theta}
\end{align*}
for some $\theta=\theta(d,q)\in(0,1)$. Since $M(2^*)\leq M(\infty)$ we may simply divide by $M(\infty)^{\theta}$ to get $M(\infty)\leq C_3\rho(q)$. Now the claim follows as in the case $2q<d$.
\hfill $\Box$

\subsection{Domination in $L_{q^{**}}(\Omega)$}\label{Domination in Lr} \textit{Proof of Theorem \ref{Dominationtheorem}.}
We only show the assertion about the Dirichlet problem (\ref{DP1}). The statement about (\ref{DP0}) is proved similarly.

Let us define the operator $\A$ with domain $\{v\in H^1_0(\Omega): \Delta v \in L_2(\Omega)\}$ which acts as $-\Delta - {\lambda_1+\varepsilon}$. By the Poincar\'e-Steklov inequality this operator is invertible.

Furthermore let us define the nonlinear but continuous operator $S:H^1_0(\Omega)\rightarrow L_2(\Omega)$ by $Sv=L\abs{\nabla v}+h$. Then (\ref{DP1}) is equivalent to the fixed point problem
\begin{align*}
	v = \A^{-1}Sv =: Tv
\end{align*}
for the operator $T:H^1_0(\Omega)\rightarrow H^1_0(\Omega)$. We apply Banach's contraction mapping principle to show existence and uniqueness of a solution. It is important to choose an appropriate norm on $H^1_0(\Omega)$ which makes $T$ a contraction mapping. 

Let $v_1,v_2\in H^1_0(\Omega)$ be arbitrary, $u_i:= Tv_i$ for $i=1,2$ and $v:=v_2-v_1$ and $u:=u_2-u_1$. Then
\begin{align*}
	\left\langle \A u, \varphi \right\rangle_{L_2(\Omega)} = \left\langle Sv_2-Sv_1, \varphi \right\rangle_{L_2(\Omega)} \quad \forall \varphi\in H^1_0(\Omega).
\end{align*}
Let $\varepsilon_1,\varepsilon_2\geq0$ with $\varepsilon=\varepsilon_1+\varepsilon_2$ for which there is a $\delta_1\in(0,1]$ such that $\varepsilon_1=\delta_1 \lambda_1$. If $\lambda_1=0$ we choose $\delta_1=1$. From the last equation with $\varphi=u\in H^1_0(\Omega)$ follows 
\begin{align*}
	\delta_1\norm{\nabla u}_{L_2(\Omega)}^2 + \varepsilon_2\norm{u}_{L_2(\Omega)}^2 \leq \alpha\delta_1\norm{\nabla v}_{L_2(\Omega)}^2 + \frac{L^2}{4\alpha\delta_1}\norm{u}_{L_2(\Omega)}^2
\end{align*}
for all $0<\alpha<1$ as in the derivation of (\ref{Bumm}) in the proof of Theorem \ref{H1}. If we set $\varepsilon_1=\min\{\varepsilon/2, \lambda_1\}$ again we see that $4\delta_1\varepsilon_2>L^2$ and thus there exists an $\alpha$, maybe close to $1$, such that $\varepsilon':=\varepsilon_2-\frac{L^2}{4\alpha\delta_1}>0$. Thus we proved that $T$ is a contraction mapping with contraction constant $\sqrt{\alpha}$ with respect to the norm $\left(\delta_1\norm{\nabla u}_{L_2(\Omega)}^2 + \varepsilon'\norm{u}_{L_2(\Omega)}^2\right)^{1/2}$. This shows that the Dirichlet problem (\ref{DP1}) has a unique solution $\vup$.

It remains to show $\vup\geq0 \text{ and } u\leq \vup$ for every solution $u$ of (\ref{DPklasse}). We only show the more difficult second assertion. Therefore let $w=(u-\vup)_+\in H^1_0(\Omega)$ and $g=(\lambda_1+\mu-\varepsilon)\vup + L\abs{\nabla\vup} + h\geq 0$. Let $w$ serve as a test function for 
\begin{align*}
  \left\{
  \begin{array}{rlr}
    -\Delta (u-\vup) + \mu (u-\vup) & \leq tb_{\sigma}(x,u,\nabla u)\chi_{\omega}(x) - g(x) & (x\in\Omega)\\
    u(x)-\vup(x) & \leq 0 & (x\in \partial \Omega).
  \end{array} 
	\right.
\end{align*}
It follows that
\begin{align*}
	\norm{\nabla w}_{L_2(\Omega)}^2 + \mu \norm{w}_{L_2(\Omega)}^2 &\leq t\int_{\omega} (b_{\sigma}(x,u_+,\nabla u_+)-g)w dx \\
	&\leq t\int_{\omega} (\lambda_1+\mu-\varepsilon)w^2+L\abs{\nabla w}w dx \\
	&\leq (\lambda_1+\mu-\varepsilon)\norm{w}_{L_2(\Omega)}^2 + L \norm{\nabla w}_{L_2(\Omega)}\norm{w}_{L_2(\Omega)}.
\end{align*}
As in the proof of Proposition \ref{uniqueness Proposition} we deduce $w=0$. This means $u\leq\vup$.
\hfill $\Box$


\section{Proof of the main theorem}\label{Proof of Theorem}
Now we prove Theorem \ref{Haupttheorem}. Observe that the conditions (\ref{coercive}) and (\ref{Wachstumsbedingung}) in conjunction with the Theorems \ref{Lr} and \ref{Dominationtheorem} imply
\begin{align}\label{Zweite Wachstumsbedingung NEU}
	\abs{f(x,v(x),\xi)} \leq f_0(x) + L_0\abs{\xi} \quad \forall \vdown\leq v\leq\vup, \xi\in\R^d,
\end{align}
for some function $f_0\in L_2(\Omega), f_0\geq0$. 

We divide the proof into three steps. In the first step we do not consider the actual Dirichlet problem, but a truncated version of it. This truncation procedure makes it possible to apply Schaefer's fixed point theorem in the locally convex space $H^1_{loc}(\Omega)$ to achieve the existence of a sequence of weak solutions of such truncated problems. This will be the second step. In the last step we show that a subsequence converges to a weak solution of (\ref{DP}). Such a pattern was applied in \cite{ArendtChill2010} to get an existence theorem for a parabolic equation on an arbitrary open set.

\textbf{(i) Truncation.}
It is well known that there is an increasing sequence of open sets $(\Omega_k)$ with $C^{\infty}$-boundary such that $\overline{\Omega_k}$ is compact and included in $\Omega_{k+1}$ and $\Omega$, and such that their union is $\Omega$. Let $\sigma:\Omega\times\R\rightarrow\R$ be as in (\ref{DPklasse}) but with $v_0=\vdown$ and $v_1=\vup$. For $\mu>0$ such that $\lambda_1+\mu-\varepsilon\geq 0$ let us define the Caratheodory function $b_{\sigma}$ as in (\ref{DPklasse}) and $b(x,s,\xi) = f(x,s,\xi) + \mu s$. For $v\in H^1_{loc}(\Omega)$ we consider the following truncated Dirichlet problem
\begin{align}
\label{truncated equation}
  \left\{
  \begin{array}{rlr}
    -\Delta u + \mu u & = b_{\sigma}(x,v,\nabla v)\chi_{\Omega_k}(x) & (x\in\Omega)\\
    u(x) & = 0 & (x\in\partial \Omega).
  \end{array} 
	\right.
\end{align} 

\textbf{(ii) Schaefer's fixed point argument.}
From (\ref{Zweite Wachstumsbedingung NEU}) we deduce that the Nemytskii operator $v\mapsto b_{\sigma}(x,v,\nabla v)\chi_{\Omega_k}$ is continuous from $H^1_{loc}(\Omega)$ to $L_2(\Omega)$ and maps $\norm{\cdot}_{H^1(\Omega_k)}$-bounded\footnote{A subset of $H^1_{loc}(\Omega)$ is called $\norm{\cdot}_{H^1(\Omega_k)}$-bounded if it is included in the ball $\{u\in H^1_{loc}(\Omega): \norm{u}_{H^1(\Omega_k)}< r\}$ for some $r>0$.} sets of $H^1_{loc}(\Omega)$ into bounded sets in $L_2(\Omega)$ (cf. \cite[Theorem 19.2]{Weinberg1956}). Thus (\ref{truncated equation}) defines an operator 
\begin{align*}
T_k:H^1_{loc}(\Omega)\rightarrow H^1_0\cap H^2_{loc}(\Omega)\hookrightarrow H^1_{loc}(\Omega), \, \text{by } T_kv=u.
\end{align*} 
It is a continuous mapping on the locally convex space $H^1_{loc}(\Omega)$ and satisfies condition (ii) in Theorem \ref{Schaefers Theorem} for $X=H^1_{loc}(\Omega)$, $\norm{\cdot}=\norm{\cdot}_{H^1(\Omega_k)}$ and every $\rho>0$. Let 
\begin{align*}
	\Schaefer_k = \{u\in H^1_{loc}(\Omega): u=tT_k u \text{ for some } 0<t\leq 1\}
\end{align*}
be the Schaefer set with respect to $T_k$. If we can show that $\Schaefer_k$ is $\norm{\cdot}_{H^1(\Omega_k)}$-bounded then we get the existence of at least one fixed point $u_k=T_ku_k$ for every $k$ by Theorem \ref{Schaefers Theorem}. By Theorem \ref{H1} an even stronger assertion is true:
\begin{align*}
	u\in\Schaefer_k \Rightarrow \norm{u}_{H^1_0(\Omega)} \leq \rho,
\end{align*}
for some constant $\rho$ which does not depend on $k$. Therefore we get the existence of an $H^1_0(\Omega)$-bounded sequence $(u_k)$ of fixed points, as desired. 

\textbf{(iii) Convergence to a solution.} 
By passing to a subsequence if necessary, we may assume that $u_k$ converges weakly to some function $u$ in $H^1_0(\Omega)$. Furthermore we may assume that this convergence is also true in the pointwise sense (almost everywhere), since the embedding $H^1_{loc}(\Omega)\hookrightarrow L_{2,loc}(\Omega)$ is compact. By Theorem \ref{Dominationtheorem} we know that $\vdown\leq u_k\leq \vup$ for all $k$ and therefore these functions satisfy
\begin{align}
\label{Finale}
  \left\{
  \begin{array}{rlr}
    -\Delta u_k + \mu u_k & = b(x,u_k,\nabla u_k)\chi_{\Omega_k}(x) & (x\in\Omega)\\
    u_k(x) & = 0 & (x\in\partial \Omega).
  \end{array} 
	\right.
\end{align} 
By (\ref{Zweite Wachstumsbedingung NEU}) the estimate $\norm{b(\cdot,u_k,\nabla u_k)}_{L_2(\Omega)}\leq \norm{f_0}_{L_2(\Omega)} + \mu\rho + L_0\rho$ is true for all $k$. Thus (\ref{Finale}) implies that $(u_k)$ is also bounded with values in the domain of the Dirichlet-Laplacian which is embedded into $H^2_{loc}(\Omega)$. Thus we may assume that $(u_k)$ also strongly converges in $H^1_{loc}(\Omega)$ to $u$. Again by \cite[Theorem 19.2]{Weinberg1956} this implies that $b(x,u_k,\nabla u_k)$ strongly converges in $L_{2,loc}(\Omega)$ to $b(x,u,\nabla u)\in L_2(\Omega)$. The arguments above allow us to take the limit $k\rightarrow\infty$ in (\ref{Finale}) which shows that $u$ is a weak solution of (\ref{DP}).
\hfill $\Box$


\section{Can the bound for $L$ be improved?} \label{Gegenbeispiel}
It is clear that in general we lose the existence of solutions for (\ref{DP}) if $\varepsilon$ in condition (\ref{coercive}) is allowed to be zero. Now we ask
\begin{itemize}
	\item Do we also lose the existence of solutions if $L\geq \Lmax$?
	\end{itemize}
Unfortunately we do not know the answer of this question. In some standard situations it is not difficult to adopt the proof of Theorem \ref{H1} (which exclusively uses that $L$ is bounded from above) to get existence of solutions. For example if $f(x,s,\xi)$ is of the form $f_1(x,s)-b(x)\cdot\xi$ for some $C^1$ vector field $b$ with $\nabla\cdot b\leq 0$, then one can prove a priori bounds in $H^1_0(\Omega)$ in the sense of Theorem \ref{H1} where $C$ does not depend on $b$. This is due to $2\int_{\Omega}b\cdot(\nabla u) u dx\geq \int_{\Omega} \nabla\cdot{(bu^2)}dx = 0$ if $u\in H^1_0(\Omega)$. We remark that the structure condition $\nabla\cdot b\leq 0$ arises in applications (see \cite{NazarovUraltseva2012} and references therein). 

The existence theorem in section \ref{Main results} is stated only for \textit{real valued} boundary value problems. The only reason why we cannot extend it to elliptic \textit{complex valued} problems (or even \textit{systems}) is that Theorem \ref{Dominationtheorem} (a priori domination) does not extend to this situation. However, if we strengthen (\ref{Gamma bei unendlich}), an existence theorem can be formulated for complex problems. Keeping this in mind we now consider a complex valued problem
\begin{align}\label{L counter example}
  \left\{
  \begin{array}{rlr}
	-\Delta u + ib\cdot \nabla u - (\lambda_1-\varepsilon)u &= g(x) & (x\in\Omega),\\
	u(x) &= 0 & (x\in \partial \Omega)
	\end{array} 
	\right.
\end{align}
where $b\in\R^d$ is a constant vector, $g\in L_2(\Omega)$ and show 

\begin{theorem}
Let $r\geq \Lmax$ and $d\geq 2$. Then there exists $b\in\R^d$ with $\abs{b}=r$ an open set $\Omega$ and $g\in L_2(\Omega)$ such that (\ref{L counter example}) has no solutions.
\end{theorem}

\textit{Proof.} It is convenient to distinguish three cases: (i) $\lambda_1=0$, (ii) $\varepsilon\leq2\lambda_1$ and (iii) $0<2\lambda_1\leq\varepsilon$.

\textit{Case (iii): $0<2\lambda_1\leq\varepsilon$.} We choose $\Omega=(0,\sqrt{\lambda_1}\pi)\times\R^{d-1}$. We write $x=(x_1,x')\in \Omega$ and similarly $b=(b_1,b')$. After taking Fourier transform
\begin{align*}
	\tilde{u}(\zeta,\xi) = \int_0^{\sqrt{\lambda_1}\pi} e^{-i x_1 \zeta} u(x_1,x') dx_1 + \int_{\R^{d-1}} e^{-i x' \xi} u(x_1,x') dx'
\end{align*}
for $\zeta\in \sqrt{\lambda_1}\Z$ and $\xi\in\R^{d-1}$ we see that (\ref{L counter example}) is equivalent to
\begin{align}\label{L counter example transformed}
	\underbrace{(\zeta^2 + \abs{\xi}^2 - b_1\zeta - b'\cdot\xi -(\lambda_1 - \varepsilon))}_{=:p(\zeta,\xi)}\tilde{u} = \tilde{g} \in l_2(\sqrt{\lambda_1}\Z)\otimes_2 L_2(\R^{d-1}),
\end{align}
where a solution $\tilde u$ must necessarily belong to $ l_2(\sqrt{\lambda_1}\Z)\otimes_2 L_2(\R^{d-1})$. It is possible to find $b\in \R^d$ such that
\begin{align}\label{Wahl von b}
	b_1 = 2\sqrt{\lambda_1} \text{ and } \abs{b} = r. 
\end{align}
Thus $b/2\in \sqrt{\lambda_1}\Z$ and $p(b/2)\leq 0$. Therefore we can find $\xi\in\R^{d-1}$ such that $p(\sqrt{\lambda_1},\xi)=0$, and it is possible to find $\tilde g$ such that (\ref{L counter example transformed}) has no solution. As a consequence (\ref{L counter example}) has no solutions for the corresponding $g$.

\textit{Case (i): $\lambda_1=0$.} When we interpret $\Omega$ from the third case as $\R^d$, then the above argumentation also works in this case if we choose $b$ with $\abs{b}=r\geq 2\sqrt{\varepsilon}$ arbitrary and note that $p(\sqrt{\varepsilon}b/\abs{b})\leq 0$.

\textit{Case (ii): $\varepsilon\leq 2\lambda_1$.} Instead of (\ref{Wahl von b}) we choose $b$ with $b_1=\varepsilon/\sqrt{\lambda_1}$ and $\abs{b}=r\geq\varepsilon/\sqrt{\lambda_1}$ arbitrary and note that $p(\sqrt{\lambda_1},0)=0$.
\hfill $\Box$

\begin{remark}
We only used the fact that a weak solution of (\ref{L counter example}) has to be in $H^1(\Omega)$ and satisfies the equation in the sense of distributions. Therefore this example works not only for Dirichlet boundary conditions.
\end{remark}


\section{Robin boundary conditions}\label{Robin BVPs}
We can generalize the problem (\ref{DP}) by considering more general boundary conditions than merely Dirichlet boundary conditions. In the following we treat Robin boundary conditions, which can be defined on arbitrary open sets (cf. \cite{ArendtWarma2003,ArendtWarma2003(2),Daners2000}). As a result of this section it turns out that a version of Theorem \ref{Haupttheorem} remains true for Robin boundary conditions under a strict positivity assumption. For general Robin boundary condition we need an additional assumption (cf. (\ref{Einbettung})) which replaces the Sobolev embedding theorem which is necessary to establish the a priori bounds in the Lebesgue spaces by Moser iteration. 

\subsection{Definition of generalized Robin boundary conditions}\label{Allgemeine Form} Let $\mu$ be a (positive) measure on the Borel $\sigma$-algebra of $\partial\Omega$. We define the positive form
\begin{align}\label{Definition Form}
	a(u,v) = \int_{\Omega} \nabla u \nabla v dx + \int_{\partial\Omega} uv d\mu
\end{align}
with domain $D(a)=\{u\in H^1(\Omega)\cap C_c(\overline{\Omega})| \int_{\partial\Omega}\abs{u}^2d\mu<\infty\}$. We write $a(u):=a(u,u)$ and equip $D(a)$ with the norm $(a(u)+\norm{u}_{L_2(\Omega)}^2)^{1/2}$. If $a$ is closable (i.e. the completion of $D(a)$ embeds injectively into $L_2(\Omega)$) we may define a self-adjoint operator by:
\begin{align*}
	D(\A) = \{u\in V| \exists g\in L_2(\Omega)\forall v \in V: \form(u,v) = \skprd{g}{v}_{L_2(\Omega)}\}, \,
	\A u = -g,
\end{align*}
where $\form$ denotes the closure of $a$. This operator acts as $\A u= \Delta u \in L_2(\Omega)$. Thus $D(\A)\subset H^2_{loc}(\Omega)$, i.e. local maximal regularity does not depend on the specific boundary conditions (see \cite[Theorem 8.8]{GilbargTrudinger2001}). Even if $a$ is not closable there exists a (unique) \textit{maximal} closable positive form $a_r$ smaller than $a$ in the following sense: The form $a_r$ is \textit{smaller} than $a$, that is $D(a)\subseteq D(a_r)$ and $a_r(u)\leq a(u)$ for every $u\in D(a)$, and every closable positive form $b$ which is smaller than $a$ is also smaller than $a_r$ (see \cite[Supplementary material, Theorem S.15]{ReedSimon_I}).

In \cite{ArendtWarma2003(2)} the authors give a characterization of $a_r$ by means of the relative capacity with respect to $\Omega$. The \textit{relative capacity} for a (not necessarily Borel-) subset $A\subseteq\overline{\Omega}$ is defined by 
\begin{align*}
	\relCap(A) = \inf\{\int_{\Omega}(\abs{\nabla u}^2+\abs{u}^2) dx\,| \, &u\in \tilde{H}^1(\Omega), \exists O\subseteq\overline{\Omega} \text{ relatively open}:\\  &A\subseteq O \text{ and } u\geq 1 \text{ a.e. on } O \}
\end{align*}
Here $\tilde{H}^1(\Omega)$ denotes the closure of $H^1(\Omega)\cap C(\overline{\Omega})$ in $H^1(\Omega)$. The relative capacity is an outer measure. A property is said to hold \textit{relatively quasi-everywhere} (r.q.e.) if it holds on $\overline{\Omega}\backslash N$ where $N$ is a set with $\relCap(N)=0$. Every $u\in \tilde{H}^1(\Omega)$ has a r.q.e. unique \textit{relatively quasi-continuous} representative which we denote by $\tilde{u}$. This is a function $\tilde{u}=u$ a.e. on $\Omega$ such that for each $\varepsilon>0$ there exists a relatively open subset $\omega\subseteq\overline{\Omega}$ with $\relCap(\overline{\Omega}\backslash\omega)\leq\varepsilon$ and $\tilde{u}$ restricted to $\omega$ is continuous (see \cite[Chapter I, Theorem 8.2.1]{BouleauHirsch1991}). Let $\Gamma_{\mu}=\{x\in\partial\Omega| \exists r>0: \mu(\{y\in\partial\Omega|\abs{x-y}<r\})<\infty\}$ be the maximal open set where $\mu$ is locally finite. For a Borel set $S\subseteq\Gamma_{\mu}$ we define $a_S$ to be the positive form with domain $D(a_S)=D(a)$ which is given by (\ref{Definition Form}) where the boundary integral over $\partial\Omega$ is replaced by the same integral over $S$. Note that $a_S$ is smaller then $a$. The Borel set $S\subseteq \Gamma_{\mu}$ is called $\mu$-\textit{admissible} iff $\relCap(A)=0$ implies $\mu(A)=0$ for each Borel-subset $A$ of $S$. 

\begin{theorem}[\cite{ArendtWarma2003,ArendtWarma2003(2)}]\label{ArWa}
There exists a $\mu$-admissible set $S\subseteq\Gamma_{\mu}$ such that $\relCap(\Gamma_{\mu}\backslash S)=0$ and $a_S=a_r$. Moreover
\begin{align*}
	D(\form_{S}) = \{u\in \tilde{H}^1(\Omega)| \tilde{u}=0 \text{ r.q.e on } \partial\Omega\backslash S \text{ and } \int_{S}\abs{\tilde{u}}^2 d\mu <\infty \} \\
	\form_S(u,v) = \int_{\Omega} \nabla u \nabla v dx + \int_S \tilde{u}\tilde{v} d\mu.
\end{align*}
$S$ is unique up to a $\mu$-null set.
\end{theorem}
In \cite[Example 4.3]{ArendtWarma2003(2)} the authors constructed a bounded domain $\Omega$ such that $\partial\Omega$ is not admissible for the $(d-1)$-dimensional Hausdorff-measure $\sigma$ although it has finite measure with respect to $\sigma$. This shows that the maximal admissible set $S$ given by the above Theorem does not in general coincide with $\Gamma_{\mu}$.

\subsection{A generalized existence theorem}
Let $\mu$ be a measure on the Borel $\sigma$-algebra of $\partial\Omega$ and $S$ the maximal $\mu$-admissible set from Theorem \ref{ArWa}. We consider the boundary value problem
\begin{align}
\label{BVP}
  \left\{
  \begin{array}{rlr}
    -\Delta u(x) & = f(x,u(x),\nabla u(x)) & (x\in\Omega)\\
     u(x) &= 0 \quad &(x\in\partial\Omega\backslash S) \\
	   \frac{\partial u}{\partial \nu}(x) + u(x)d\mu(x) &= 0 \quad &(x\in S).
  \end{array} 
	\right.
\end{align}
We call $u$ a \textit{weak solution} of this problem iff $u\in D(\form_S)$ and $f(x,u,\nabla u)\in L_2(\Omega)$ such that
\begin{align*}
	\form_S(u,\varphi) = \int_{\Omega} f(x,u,\nabla u)\varphi dx \quad \forall \varphi \in D(\form_S).
\end{align*}
Let $2\leq \hat{d}<\infty$ be a real number. We assume that
\begin{align}
	\label{Einbettung}
	D(\form_S) \hookrightarrow L_p(\Omega) 	
		\left\{ 
	\begin{array}{r}
	\text{for } p=2\hat{d}/(\hat{d}-2)_+ \text{ if } \hat{d}\neq 2 \\ 
	\text{ and for all } p< \infty \text{ else.}
	\end{array}
	\right.
\end{align}

All conditions on $f$ remain unchanged as compared to Section \ref{Introduction} but (\ref{Gamma bei unendlich}) has to be satisfied if and only if $q< \hat{d}/2$, where we redefine $q^{**}=q\hat{d}/(\hat{d}-2q)_+$. The (for simplicity) excluded value for $q$ is $\hat{d}/2$ instead of $d/2$. Furthermore, $\lambda_1\in[0,\infty)$ is the optimal constant such that $\lambda_1\norm{u}^2_{L_2(\Omega)}\leq \form_S(u)$ is true for all $u\in D(\form_S)$. 
\begin{theorem}\label{Haupttheorem mit Randbedingungen}
Under the above conditions the boundary value problem (\ref{BVP}) has at least one weak solution.
\end{theorem}
It is not difficult to generalize Theorems \ref{H1}, \ref{Lr} and \ref{Dominationtheorem} and the argumentation in Section \ref{Proof of Theorem} to this more general setting. In fact, (\ref{Einbettung}) replaces the Sobolev-embedding $H^1_0(\Omega)\hookrightarrow L_{2^{*}}(\Omega)$ which has to be used in step (i) of the proof of Theorem \ref{Lr} and the following lemma guarantees that the test function argument described there also works in our more general situation: 

\begin{lemma}\label{b Lemma}
Let $\beta\geq 1$ and $t_0>0$. Let $H(s)=\abs{s}^{\beta-1}s$ for all real $\abs{s}\leq t_0$. Extend $H$ affine linearly to a $C^1$-function. Then
\begin{align*}
	u\in D(\form_S) \Longrightarrow H(u) \in D(\form_S) \text{ and } H(\tilde{u}) = H(u)\,\tilde{ } 
\end{align*}
where $\tilde{u}$ and $H(u)\,\tilde{ }$ denote the relatively quasi-continuous representatives of $u$ and $H(u)$ respectively.
\end{lemma}
\textit{Proof.}
Essentially we only have to proof that $H(\tilde{u})$ is the relatively quasi-continuous representative of $H(u)$. But this is easy since if $\tilde{u}$ is continuous on some set then so is $H(\tilde{u})$.
\hfill $\Box$
\\
\\
It is only left to find examples where (\ref{Einbettung}) is true for some $\hat{d}$. Of course, if $\partial\Omega$ satisfies a uniform lipschitz condition then it is true for $\hat{d}=d$ by the Sobolev embedding theorem. But there are also other situations where (\ref{Einbettung}) is true. For this purpose set $d\mu=\beta(x)d\sigma$ where $\sigma$ is the $d-1$-dimensional Hausdorff-measure restricted to $\partial\Omega$. The Borel function $\beta$ is bounded from below, that is $\beta(x)\geq\beta_0>0$ for some constant $\beta_0$. Then it is a consequence of an inequality due to Maz'ya that (\ref{Einbettung}) holds with $\hat{d}=2d$. We refer to \cite[Chapter 5]{ArendtWarma2003(2)} for the short proof.

\begin{remark} One could ask if the $L_{q^{**}}$ a priori bounds derived from (\ref{Einbettung}) by Moser iteration are the best possible. Indeed they are: \cite[Theorem 5.11.]{Daners2000}. 
\end{remark}

\subsection*{Acknowledgments} I am most grateful to Ralph Chill and Alexander I. Nazarov who helped me with their advice and knowledge. I am also grateful to Arina A. Arkhipova for useful comments improving the final version of the article. My ten month internship at Saint Petersburg Department of V.A. Steklov Institute was supported by the DAAD (German Academic Exchange Service).


\bibliographystyle{alpha}
\bibliography{Referenzen}

\end{document}